\documentclass[reqno,a4paper]{amsart}
\usepackage{url,amssymb}
\usepackage[hypertex]{hyperref}
\allowdisplaybreaks[4]
\theoremstyle{plain}
\newtheorem{thm}{Theorem}
\theoremstyle{definition}
\newtheorem{dfn}{Definition}
\theoremstyle{remark}

\DeclareMathOperator{\td}{d\mspace{-2mu}}

\begin{document}

\title[Logarithmically absolutely monotonic functions]
{A property of logarithmically absolutely monotonic functions and the logarithmically complete monotonicity of a power-exponential function}

\author[F. Qi]{Feng Qi}
\address[F. Qi]{Research Institute of Mathematical Inequality Theory, Henan Polytechnic University, Jiaozuo City, Henan Province, 454010, China}
\email{\href{mailto: F. Qi <qifeng618@gmail.com>}{qifeng618@gmail.com},
\href{mailto: F. Qi <qifeng618@hotmail.com>}{qifeng618@hotmail.com},
\href{mailto: F. Qi <qifeng618@qq.com>}{qifeng618@qq.com}}
\urladdr{\url{http://qifeng618.spaces.live.com}}

\author[B.-N. Guo]{Bai-Ni Guo}
\address[B.-N. Guo]{School of Mathematics and Informatics, Henan Polytechnic University, Jiaozuo City, Henan Province, 454010, China}
\email{\href{mailto: B.-N. Guo <bai.ni.guo@gmail.com>}{bai.ni.guo@gmail.com}, \href{mailto: B.-N. Guo <bai.ni.guo@hotmail.com>}{bai.ni.guo@hotmail.com}}

\begin{abstract}
In the article, a notion ``logarithmically absolutely monotonic function'' is
introduced, an inclusion that a logarithmically absolutely monotonic function
is also absolutely monotonic is revealed, the logarithmically complete
monotonicity and the logarithmically absolute monotonicity of the function
$\bigl(1+\frac{\alpha}x\bigr) ^{x+\beta}$ are proved, where $\alpha$ and
$\beta$ are given real parameters, a new proof for the inclusion that a
logarithmically completely monotonic function is also completely monotonic is
given, and an open problem is posed.
\end{abstract}



\thanks{The authors were supported partially by NSF of Henan University, China}

\thanks{This paper was typeset using \AmS-\LaTeX}

\maketitle

\section{Introduction}

Recall \cite{clark-ismail.tex, clark-ismail.tex-rgmia, royal-98,haer,widder} that a function $f$ is
said to be completely monotonic on an interval $I$ if $f$ has derivative of
all orders on $I$ such that
\begin{equation}
(-1)^{k}f^{(k)}(x)\ge0
\end{equation}
for $x\in I$ and $k\ge0$.
For our own convenience, the set of the completely monotonic functions on $I$
is denoted by $\mathcal{C}[I]$.
\par
Recall also \cite{clark-ismail.tex, clark-ismail.tex-rgmia, royal-98,b^x-a^x, haer,widder} that a
function $f$ is said to be absolutely monotonic on an interval $I$ if it has
derivatives of all orders and
\begin{equation}
f^{(k-1)}(t)\ge0
\end{equation}
for $t\in I$ and $k\in\mathbb{N}$, where $\mathbb{N}$ denotes the
set of all positive integers. The set of the absolutely monotonic
functions on $I$ is denoted by $\mathcal{A}[I]$.
\par
Recall again \cite{Atanassov, clark-ismail.tex, clark-ismail.tex-rgmia, compmon2, minus-one,
auscmrgmia, e-gam-rat-comp-mon} that a positive function $f$ is said to be
logarithmically completely monotonic on an interval $I$ if its logarithm $\ln
f$ satisfies
\begin{equation}
(-1)^k[\ln f(x)]^{(k)}\ge0
\end{equation}
for $k\in\mathbb{N}$ on $I$. Similar
to above, the set of the logarithmically completely monotonic functions on $I$
is denoted by $\mathcal{C_L}[I]$.
\par
The famous Bernstein-Widder's Theorem \cite[p.~161]{widder} states that
$f\in\mathcal{C}[(0,\infty)]$ if and only if there exists a bounded and
nondecreasing function $\mu(t)$ such that
\begin{equation}\label{berstein-1}
f(x)=\int_0^{\infty}{e}^{-xt}\td \mu(t)
\end{equation}
converges for $0<x < \infty$, and that $f(x)\in\mathcal{A}[(0,\infty)]$ if and
only if there exists a bounded and nondecreasing function $\sigma(t)$ such
that
\begin{equation}
f(x)=\int_0^{\infty}{e}^{xt}\td \sigma(t)
\end{equation}
converges for $0\le x <\infty$.
\par
In \cite{CBerg,compmon2,minus-one,auscmrgmia,e-gam-rat-comp-mon,haer} and many
other references, the inclusions $\mathcal{C_L}[I]\subset\mathcal{C}[I]$ and
$\mathcal{S}\subset\mathcal{C_L}[(0,\infty)]$ were revealed implicitly or
explicitly, where $\mathcal{S}$ denotes the class of Stieltjes transforms.
\par
The class $\mathcal{C_L}[(0,\infty)]$ is characterized in
\cite[Theorem~1.1]{CBerg} explicitly and in \cite[Theorem~4.4]{horn}
implicitly: $f\in\mathcal{C_L}[(0,\infty)]\Longleftrightarrow
f^\alpha\in\mathcal{C}$ for all $\alpha>0$ $\Longleftrightarrow
\sqrt[n]{f}\,\in\mathcal{C}$ for all $n\in\mathbb{N}$. In other words,  the
functions in $\mathcal{C_L}[(0,\infty)]$ are those completely monotonic
functions for which the representing measure $\mu$ in \eqref{berstein-1} is
infinitely divisible in the convolution sense: For each $n\in\mathbb{N}$ there
exists a positive measure $\nu$ on $[0,\infty)$ with $n$-th convolution power
equal to $\mu$.
\par
To the best of our knowledge, the terminology ``logarithmically completely
monotonic function'' and some properties of it appeared firstly without
explicit definition in \cite{Atanassov}, re-coined independently with explicit
definition by the first author in \cite{minus-one}, and published formally in
\cite{compmon2}. Since then, a further deep investigation on the
logarithmically completely monotonic functions was explicitly carried out in
\cite{CBerg} and a citation of the logarithmically completely monotonic
functions appeared in \cite{grin-ismail}.
\par
It is said in \cite{CBerg} that ``In various papers complete monotonicity for
special functions has been established by proving the stronger statement that
the function is a Stieltjes transform''. It is also said in
\cite{berg-pedersen-jcam} that ``In concrete cases it is often easier to
establish that a function is a Stieltjes transform than to verify complete
monotonicity''. Because the logarithmically completely monotonic functions must
be completely monotonic, in order to show some functions, especially the
power-exponential functions or the exponential functions, are completely
monotonic, maybe it is sufficient and much simpler to prove their
logarithmically complete monotonicity or to show that they are Stieltjes
transforms, if possible. These techniques have been used in
\cite{alzer-mc-97,alzer1,Alzer1,anderson-qiu, Bustoz-and-Ismail, chen-qi-log-jmaa,ajmaa,global, global-rgmia, clark-ismail-orig, elber-laforgia-pams, grin-ismail,Ismail-Lorch-Muldoon,li-zhao-chen-beog, a33, sandor-gamma-3.tex, sandor-gamma-3.tex-rgmia, clark-ismail.tex, clark-ismail.tex-rgmia, traffic, traffic-jkms, cheung-qi-spec-func, auscmrgmia, e-gam-rat-comp-mon, Qi-Li-two, Niu-Qi-gamma-jmaat.tex, Niu-Qi-gamma-jmaat.tex-rgmia, qi-yang-li, haer, zhang-su-ling} and many other articles. It can be imagined that, if there is no the inclusion relationships between the sets of completely monotonic functions, logarithmically completely monotonic functions and Stieltjes transforms, it would be very complex, difficult, even impossible, to verify some power-exponential functions to be completely monotonic.
\par
It is worthwhile to point out that, in most related papers before, although
the logarithmically completely monotonicities of some functions had been
established essentially, but they were stated using the notion ``completely
monotonic function'' instead of ``logarithmically completely monotonic
function''. Because a completely monotonic function may be not logarithmically
completely monotonic, in my opinion, it should be better that many results or
conclusions of completely monotonic functions are rewritten or restated using
the term ``logarithmically completely monotonic function''.
\par
The main results of this paper are as follows.
\par
Similar to the definition of the logarithmically completely monotonic
function, we would like to coin a notion ``logarithmically absolutely
monotonic functions''.

\begin{dfn}
A positive function $f$ is said to be logarithmically absolutely monotonic on
an interval $I$ if it has derivatives of all orders and $[\ln f(t)]^{(k)}\ge
0$ for $t\in I$ and $k\in\mathbb{N}$.
\end{dfn}

For our own convenience, the set of the logarithmically absolutely monotonic
functions on an interval $I$ is denoted by $\mathcal{A_L}[I]$.
\par
Similar to the inclusion $\mathcal{C_L}[I]\subset\mathcal{C}[I]$ mentioned
above, the logarithmically absolutely monotonic functions have the following
nontrivial property.

\begin{thm}\label{thm-abs-log}
A logarithmically absolutely monotonic function on an interval $I$ is also
absolutely monotonic on $I$, but not conversely. Equivalently,
$\mathcal{A_L}[I]\subset\mathcal{A}[I]$ and
$\mathcal{A}[I]\setminus\mathcal{A_L}[I]\ne\emptyset$.
\end{thm}

This theorem hints us that, in order to show some functions, especially the
power-exponential functions or the exponential functions, are absolutely
monotonic, maybe it is much simpler or easier to prove the stronger statement
that they are logarithmically absolutely monotonic.
\par
Let
\begin{equation}
F_{\alpha,\beta}(x)=\biggl(1+\frac{\alpha}x\biggr)^{x+\beta}
\end{equation}
for $\alpha\ne0$ and either $x>\max\{0,-\alpha\}$ or $x<\min\{0,-\alpha\}$. In
\cite{guobn1, guobn2, guo-qi-k-ye, sxtb, amm, qif, miqschur, lim-mei-tan,
AMEN,schur-complete, xuxq} (See also the related contents in
\cite{kuang-2rd,kuang-3rd}), the sufficient and necessary conditions such that
the function $F_{\alpha,\beta}(x)$, its simplified forms, its variants and
their corresponding sequences are monotonic are obtained.
\par
In \cite[Theorem~1.2]{Niu-Qi-gamma-jmaat.tex}, \cite[Theorem~2]{Niu-Qi-gamma-jmaat.tex-rgmia} and \cite{zhang-su-ling}, it was
proved that $F_{\alpha,\beta}(x)\in\mathcal{C_L}[(0,\infty)]$ for $\alpha>0$
and $\beta\in\mathbb{R}$ if and only if $2\beta\ge\alpha>0$. From
$\mathcal{C_L}[I]\subset\mathcal{C}[I]$ it is deduced that the function
\begin{equation}
\biggl(1+\frac{\alpha}x\biggr)^{x+\beta}-e^\alpha\in\mathcal{C}[(0,\infty)]
\end{equation}
if and only
if $0<\alpha\le2\beta$, which is a conclusion obtained in \cite{alzer-berg}.
\par
Now it is natural to pose a problem: How about the logarithmically complete or
absolute monotonicity of the function $F_{\alpha,\beta}(x)$ for all real
numbers $\alpha\ne0$ and $\beta$ in the interval $(-\infty,\min\{0,-\alpha\})$
or $(\max\{0,-\alpha\},\infty)$? The following Theorem~\ref{thm4} answers this
problem.

\begin{thm}\label{thm4}
\begin{enumerate}
\item
For $\alpha<0$, $F_{\alpha,\beta}(x) \in\mathcal{C_L} [(-\alpha,\infty)]$ if
and only if $\beta\le\alpha$ and $\frac1{F_{\alpha,\beta}(x)}\in\mathcal{C_L}
[(-\alpha,\infty)]$ if and only if $2\beta\ge\alpha$.
\item
For $\alpha>0$, $F_{\alpha,\beta}(x)\in\mathcal{C_L}[(0,\infty)]$ if and only
if $2\beta\ge\alpha$ and $\frac1{F_{\alpha,\beta}(x)}\in\mathcal{C_L}
[(0,\infty)]$ if and only if $\beta\le0$.
\item
For $\alpha<0$, $F_{\alpha,\beta}(x)\in\mathcal{A_L}[(-\infty,0)]$ if and only
if $\beta\ge0$ and $\frac1{F_{\alpha,\beta}(x)}\in\mathcal{A_L}[(-\infty,0)]$
if and only if $2\beta\le\alpha$.
\item
For $\alpha>0$, $F_{\alpha,\beta}(x)\in\mathcal{A_L}[(-\infty,-\alpha)]$ if
and only if $2\beta\le\alpha$ and $\frac1{F_{\alpha,\beta}(x)}\in\mathcal{A_L}
[(-\infty,-\alpha)]$ if and only if $\beta\ge\alpha$.
\end{enumerate}
\end{thm}

As an immediate consequence of combining Theorem~\ref{thm4} with the inclusion
$\mathcal{C_L}[I]\subset\mathcal{C}[I]$, the following complete monotonicity
relating to the function $F_{\alpha,\beta}(x)$, which extends the
corresponding results in \cite{alzer-berg, Niu-Qi-gamma-jmaat.tex, Niu-Qi-gamma-jmaat.tex-rgmia, zhang-su-ling}, can be obtained easily.

\begin{thm}\label{abs-thm3}
\begin{enumerate}
\item
For $\alpha>0$,
\begin{equation}
\biggl(1+\frac{\alpha}x\biggr)^{x+\beta}-e^\alpha
\in\mathcal{C}[(0,\infty)]
\end{equation}
if and only if $\alpha\le2\beta$ and
\begin{equation}
\biggl(1+\frac{\alpha}x\biggr)^{-(x+\beta)}-e^{-\alpha}\in\mathcal{C}[(0,\infty)]
\end{equation}
if and only if $\beta\le0$.
\item
For $\alpha<0$,
\begin{equation}
\biggl(1+\frac{\alpha}x\biggr)^{x+\beta}-e^\alpha\in\mathcal{C}[(-\alpha,\infty)]
\end{equation}
if and only if $\beta\le\alpha$ and
$\frac1{F_{\alpha,\beta}(x)} \in\mathcal{C} [(-\alpha,\infty)]$ if and only if
$2\beta\ge\alpha$.
\end{enumerate}
\end{thm}

In \cite{CBerg} and \cite{compmon2,minus-one}, two different proofs
for the inclusion $\mathcal{C_L}[I]\subset\mathcal{C}[I]$ were
given. Now we would like to present a new proof for this inclusion
by using Fa\'a di Bruno's formula \cite{FaadiBrunos-orig,
FaadiBrunos-survey, FaadiBrunos}.

\begin{thm}[\cite{CBerg,compmon2,minus-one}]\label{cl-in-c}
A logarithmically completely monotonic function on an interval $I$ is also
completely monotonic on $I$, but not conversely. Equivalently,
$\mathcal{C_L}[I]\subset\mathcal{C}[I]$ and
$\mathcal{C}[I]\setminus\mathcal{C_L}[I]\ne\emptyset$.
\end{thm}

Now we are in a position to pose an open problem: How
to characterize the logarithmically absolutely monotonic functions as we
characterize the (logarithmically) completely monotonic functions and the
absolutely monotonic functions by Bernstein-Widder's Theorem?

\section{Proofs of theorems}

\begin{proof}[Proof of Theorem \ref{thm-abs-log}]
The Fa\'a di Bruno's formula \cite{FaadiBrunos-orig,
FaadiBrunos-survey, FaadiBrunos} gives an explicit formula for the
$n$-th derivative of the composition $g(h(t))$: If $g(t)$ and $h(t)$
are functions for which all the necessary derivatives are defined,
then
\begin{equation}\label{faa-bruno}
\frac{\td{}^n}{\td x^n}[g(h(x))]
=\sum_{\substack{1\leqslant i\leqslant n,i_k\geqslant0\\
\sum_{k=1}^n i_k=i\\ \sum_{k=1}^n ki_k=n}}\frac{n!}{\prod\limits_{k=1}^n
i_k!}g^{(i)}(h(x))\prod_{k=1}^n\biggl[\frac{h^{(k)}(x)}{k!}\biggr]^{i_k}.
\end{equation}
Applying \eqref{faa-bruno} to $g(x)=e^x$ and $h(x)=\ln f(x)$ leads to
\begin{equation}\label{abs-4}
f^{(n)}(x)=\bigl[e^{\ln f(x)}\bigr]^{(n)}=n!f(x)
\sum_{\substack{1\leqslant i\leqslant n,i_k\geqslant0\\
\sum_{k=1}^n i_k=i\\ \sum_{k=1}^n ki_k=n}}\prod_{k=1}^n\frac{\bigl\{{[\ln
f(x)]^{(k)}}\bigr\}^{i_k}}{[i_k!(k!)^{i_k}]}
\end{equation}
for $n\in\mathbb{N}$. If $f(x)\in\mathcal{A_L}[I]$, then $[\ln
f(x)]^{(k)}\ge0$ for $k\in\mathbb{N}$, and then $f^{(n)}(x)\ge0$ for
$n\in\mathbb{N}$, that means $f(x)\in\mathcal{A}[I]$.
\par
Conversely, it is clear that $0\in\mathcal{A}[I]$, but
$0\not\in\mathcal{A_L}[I]$. Therefore
$\mathcal{A}[I]\setminus\mathcal{A_L}[I]\ne\emptyset$. The proof of
Theorem~\ref{thm-abs-log} is complete.
\end{proof}

\begin{proof}[Proof of Theorem \ref{thm4}]
Calculating directly yields
\begin{align}
\label{js1} \ln F_{\alpha,\beta}(x)&=
(x+\beta)\ln\left(1+\frac{\alpha}{x}\right),\\\label{js2} [\ln
F_{\alpha,\beta}(x)]'&= \ln\left(1+\frac{\alpha}{x}\right)
-\frac{\alpha(x+\beta)}{x(x+\alpha)},\\
\label{js3} [\ln F_{\alpha,\beta}(x)]'' &=
\frac{\alpha[(2\beta-\alpha)x+\alpha\beta]}{x^2(x+\alpha)^2},
\end{align}
and
\begin{equation}\label{js4}
\lim_{x\to\pm\infty}[\ln F_{\alpha,\beta}(x)]''= 0, \qquad
\lim_{x\to\pm\infty}[\ln F_{\alpha,\beta}(x)]'=0.
\end{equation}
For $\alpha<0$ and $x>-\alpha$, in virtue of formula
\begin{equation}\label{x^r}
\frac1{x^r}=\frac1{\Gamma(r)} \int_0^\infty t^{r-1}e^{-xt}\td t
\end{equation}
for $x>0$ and $r>0$, equation \eqref{js3} can be rewritten as
\begin{equation}
\begin{split}\label{f''}
[\ln F_{\alpha,\beta}(x)]''&=\frac1{x+\alpha}-\frac1x-\frac{\beta-\alpha}
{(x+\alpha)^2}+\frac{\beta}{x^2}\\
&\triangleq\int_0^\infty[\beta-q_\alpha(t)] t\bigl(e^{\alpha
t}-1\bigr)e^{-(x+\alpha)t}\td t,
\end{split}
\end{equation}
where
\begin{equation}\label{pqu}
q_\alpha(t)=\frac{e^{\alpha t}-\alpha t-1} {t\bigl(e^{\alpha t}-1\bigr)}
=\frac{\alpha(e^{u}-u-1)} {u(e^{u}-1)}\triangleq\alpha q(u)
\end{equation}
for $t>0$ and $u=\alpha t<0$. Since $q(u)$ is decreasing in $(-\infty,0)$ with
$\lim_{u\to0^-}q(u)=\frac12$ and $\lim_{u\to-\infty}q(u)=1$, then when
$\beta\le\alpha$ the function $(-1)^i[\ln F_{\alpha,\beta}(x)]^{(i+2)}\ge0$
and when $2\beta\ge\alpha$ the function $(-1)^i[\ln
F_{\alpha,\beta}(x)]^{(i+2)}\le0$ in $(-\alpha,\infty)$ for $i\ge0$. Since
$[\ln F_{\alpha,\beta}(x)]'$ increases for $\beta\le\alpha$ and decreases for
$2\beta\ge\alpha$, considering one of the limits in \eqref{js4} shows that
$[\ln F_{\alpha,\beta}(x)]'\le0$ for $\beta\le\alpha$ and $[\ln
F_{\alpha,\beta}(x)]'\ge0$ for $2\beta\ge\alpha$. In conclusion, $(-1)^k[\ln
F_{\alpha,\beta}(x)]^{(k)}\ge0$ for $\beta\le\alpha$ and $(-1)^k[\ln
F_{\alpha,\beta}(x)]^{(k)}\le0$ for $2\beta\ge\alpha$ and $k\in\mathbb{N}$.
This means that $F_{\alpha,\beta}(x)\in\mathcal{C_L}[(-\alpha,\infty)]$ for
$\beta\le\alpha<0$ and $\frac1{F_{\alpha,\beta}(x)}\in \mathcal{C_L}
[(-\alpha,\infty)]$ for both $2\beta\ge\alpha$ and $\alpha<0$. Conversely, if
$F_{\alpha,\beta}(x)\in\mathcal{C_L}[(-\alpha,\infty)]$ for $\alpha<0$, then
$[\ln F_{\alpha,\beta}(x)]'\le0$ which can be rearranged as
\begin{equation}
\beta\le x\left[\left(1+\frac{x}\alpha\right)
\ln\left(1+\frac\alpha{x}\right)-1\right]\triangleq \theta_\alpha(x)
\end{equation}
and $\lim_{x\to(-\alpha)^+}\theta_\alpha(x)=\alpha$, thus $\beta\le\alpha$. If
$\frac1{F_{\alpha,\beta}(x)}\in\mathcal{C_L}[(-\alpha,\infty)]$ for
$\alpha<0$, then $[\ln F_{\alpha,\beta}(x)]'\ge0$ which can be rearranged as
$\beta\ge\theta_\alpha(x)\to\frac\alpha2$ as $x\to\infty$, hence
$2\beta\ge\alpha$ holds.
\par
If $\alpha>0$, the formulas \eqref{f''} and \eqref{pqu} are valid for $x>0$
and $u>0$. Since $q(u)$ is decreasing in $(0,\infty)$ with
$\lim_{u\to0^+}q(u)=\frac12$ and $\lim_{u\to\infty}q(u)=0$, by the same
argument as above, it follows easily that $F_{\alpha,\beta}(x)
\in\mathcal{C_L} [(0,\infty)]$ for $2\beta\ge\alpha$ and
$\frac1{F_{\alpha,\beta}(x)} \in \mathcal{C_L} [(0,\infty)]$ for $\beta\le0$.
Conversely, if $F_{\alpha,\beta}(x) \in\mathcal{C_L} [(0,\infty)]$ for
$\alpha>0$, then $[\ln F_{\alpha,\beta}(x)]'\le0$ which can be rewritten as
$\beta\ge\theta_\alpha(x)\to\frac\alpha2$ as $x$ tends to $\infty$; if
$\frac1{F_{\alpha,\beta}(x)}\in\mathcal{C_L} [(0,\infty)]$ for $\alpha>0$,
then $[\ln F_{\alpha,\beta}(x)]'\ge0$ which can be rewritten as
$\beta\le\theta_\alpha(x)\to0$ as $x\to0$.
\par
For $\alpha<0$ and $x<0$, it is easy to obtain
\begin{equation}
\begin{split}\label{f'''}
[\ln F_{\alpha,\beta}(x)]''&=-\frac1{-(x+\alpha)}+\frac1{-x}
-\frac{\beta-\alpha}{[-(x+\alpha)]^2}+\frac{\beta}{(-x)^2}\\
&\triangleq\int_0^\infty[\beta+p_\alpha(t)] t\bigl(1-e^{\alpha
t}\bigr)e^{xt}\td t,
\end{split}
\end{equation}
where
\begin{equation}\label{pqu'}
p_\alpha(t)=\frac{1+(\alpha t-1)e^{\alpha t}}{t(1-e^{\alpha t})}
=\frac{\alpha[1+(u-1)e^{u}]}{u(1-e^{u})}\triangleq\alpha p(u)
\end{equation}
for $t>0$ and $u=\alpha t<0$ and $p(u)$ is decreasing in $(-\infty,0)$ with
$\lim_{u\to-\infty}p(u)=0$ and $\lim_{u\to0^-}p(u)=-\frac12$. Accordingly, for
$i\ge0$ and in $(-\infty,0)$, if $\beta-\frac\alpha2\le0$ then $[\ln
F_{\alpha,\beta}(x)]^{(i+2)}\le0$, if $\beta\ge0$ then $[\ln
F_{\alpha,\beta}(x)]^{(i+2)}\ge0$. By virtue of \eqref{js4}, it is deduced
immediately that $[\ln F_{\alpha,\beta}(x)]^{(k)}\le0$ for $2\beta\le\alpha$
and $[\ln F_{\alpha,\beta}(x)]^{(k)}\ge0$ for $\beta\ge0$ and $k\in\mathbb{N}$
in $(-\infty,0)$. Conversely, if $F_{\alpha,\beta}(x)$ is logarithmically
absolutely monotonic in $(-\infty,0)$, then $[\ln F_{\alpha,\beta}(x)]'\ge0$
which can be rewritten as $\beta\ge\theta_\alpha(x)$ for $x\in(-\infty,0)$.
From $\lim_{x\to0^-}\theta_\alpha(x)=0$, it follows that $\beta\ge0$; if
$\frac1{F_{\alpha,\beta}(x)}$ is logarithmically absolutely monotonic in
$(-\infty,0)$, then $[\ln F_{\alpha,\beta}(x)]'\le0$ which can be rearranged
as $\beta\le\theta_\alpha(x)$ for $x\in(-\infty,0)$. From
$\lim_{x\to-\infty}\theta_\alpha(x)=\frac\alpha2$, it concludes that
$2\beta\le\alpha$.
\par
For $\alpha>0$ and $x<-\alpha$, the formulas \eqref{f'''} and \eqref{pqu'}
hold for $x\in(-\infty,-\alpha)$ and $u>0$. The function $p(u)$ is negative
and decreasing in $(0,\infty)$ with $\lim_{u\to0^+}p(u)=-\frac12$ and
$\lim_{u\to\infty}p(u)=-1$. Consequently, if $\beta-\frac12\alpha\le0$ then
$[\ln F_{\alpha,\beta}(x)]^{(i+2)}\ge0$ for $i\ge0$ in $(-\infty,-\alpha)$, if
$\beta-\alpha\ge0$ then $[\ln F_{\alpha,\beta}(x)]^{(i+2)}\le0$ for $i\ge0$ in
$(-\infty,-\alpha)$. In virtue of \eqref{js4}, it is concluded readily that
$[\ln F_{\alpha,\beta}(x)]^{(k)}\ge0$ for $2\beta\le\alpha$ and $[\ln
F_{\alpha,\beta}(x)]^{(k)}\le0$ for $\beta\ge\alpha$ and $k\in\mathbb{N}$ in
$(-\infty,-\alpha)$. Conversely, if $F_{\alpha,\beta}(x)$ is logarithmically
absolutely monotonic in $(-\infty,-\alpha)$, then $[\ln
F_{\alpha,\beta}(x)]'\ge0$ which can be rewritten as
$\beta\le\theta_\alpha(x)$ for $x\in(-\infty,-\alpha)$. From the fact that
$\lim_{x\to-\infty}\theta_\alpha(x)=\frac\alpha2$, it follows that
$2\beta\le\alpha$; if $\frac1{F_{\alpha,\beta}(x)}$ is logarithmically
absolutely monotonic in $(-\infty,-\alpha)$, then $[\ln
F_{\alpha,\beta}(x)]'\le0$ which can be rearranged as
$\beta\ge\theta_\alpha(x)$ for $x\in(-\infty,-\alpha)$. From the fact that
$\lim_{x\to(-\alpha)^-}\theta_\alpha(x)=\alpha$, it concludes that
$\beta\ge\alpha$. The proof of Theorem~\ref{thm4} is complete.
\end{proof}

\begin{proof}[Proof of Theorem \ref{abs-thm3}]
This follows from taking limits by L'H\^ospital rule, considering the
inclusion $\mathcal{C_L}[I]\subset\mathcal{C}[I]$ and using
Theorem~\ref{thm4}.
\end{proof}

\begin{proof}[Proof of Theorem \ref{cl-in-c}]
If $f(x)\in\mathcal{C_L}[I]$, then $(-1)^k[\ln f(x)]^{(k)}\ge0$ for
$k\in\mathbb{N}$ and
\begin{equation}\label{abs-5}
\begin{split}
\prod_{k=1}^n\bigl\{{[\ln f(x)]^{(k)}}\bigr\}^{i_k}
&=\prod_{k=1}^n(-1)^{ki_k}\bigl\{{(-1)^k[\ln f(x)]^{(k)}}\bigr\}^{i_k}\\
&=\prod_{k=1}^n(-1)^{ki_k}\prod_{k=1}^n\bigl\{{(-1)^k[\ln
f(x)]^{(k)}}\bigr\}^{i_k}\\
&=(-1)^{\sum_{k=1}^nki_k}\prod_{k=1}^n\bigl\{{(-1)^k[\ln
f(x)]^{(k)}}\bigr\}^{i_k}
\end{split}
\end{equation}
for $n\in\mathbb{N}$. Substituting \eqref{abs-5} into \eqref{abs-4} yields
\begin{equation}\label{abs-44}
(-1)^nf^{(n)}(x)=n!f(x)
\sum_{\substack{1\leqslant i\leqslant n,i_k\geqslant0\\
\sum_{k=1}^n i_k=i\\ \sum_{k=1}^n
ki_k=n}}\prod_{k=1}^n\frac{\bigl\{{(-1)^k[\ln
f(x)]^{(k)}}\bigr\}^{i_k}}{[i_k!(k!)^{i_k}]}\ge0
\end{equation}
for $n\in\mathbb{N}$. This means that $f(x)\in\mathcal{C}[I]$.
\par
Conversely, it is clear that $0\in\mathcal{C}[I]$ for any interval $I$, but
$0\not\in\mathcal{C_L}[I]$. Therefore
$\mathcal{C}[I]\setminus\mathcal{C_L}[I]\ne\emptyset$. The proof of Theorem
\ref{cl-in-c} is complete.
\end{proof}

\end{document}